\documentclass[a4paper,11pt]{amsart}
\usepackage[cp1251]{inputenc}
\usepackage[english]{babel}
\usepackage{amsmath,amssymb,euscript,amsthm,amsfonts}
\newtheorem{lemma}{Lemma}
\newtheorem{theorem}{Theorem}
\theoremstyle{definition}
\newtheorem{definition}{Definition}
\theoremstyle{remark}
\newtheorem{remark}{Remark}

\newtheorem{corollary}{Corollary}

\begin{document}

\renewcommand{\proofname}{Proof}
\makeatletter \headsep 10 mm \footskip 10 mm
\renewcommand{\@evenhead}%
{\vbox{\hbox to\textwidth{\strut \centerline{{\it Nikolaj
Glazunov}}} \hrule}}

\begin{title}
{On A.V. Malyshev's approach to Minkowski's conjecture concerning the critical determinant
 of the region $|x|^p + |y|^p < 1$  for $p > 1$}
\end{title}
\begin{author}
{Nikolaj  Glazunov}
\end{author}

\begin{abstract}
{We present  A.V. Malyshev`s approach to Minkowski`s conjecture (in Davis`s  amendment) concerning the critical 
determinant of the region $|x|^p + |y|^p <1$  for $p > 1$ and   Malyshev`s method.
In the sequel of this article we use these approach and method to present the main result.}
\end{abstract}

\maketitle

\section{Introduction}

Let
 
$$ |\alpha x + \beta y|^p + |\gamma x + \delta y|^p \leq c
 |\det(\alpha \delta - \beta \gamma)|^{p/2}, $$ 

be a diophantine inequality defined for a given real $ p >1 $;
hear $\alpha, \beta, \gamma, \delta$ are real numbers with 
$ \alpha \delta - \beta \gamma \neq 0 .$ 

H. Minkowski in his monograph~\cite{Mi:DA} raise the question
about minimum constant $c$ 
such that the inequality has integer 
solution other than origin. 
Minkowski with the help of his theorem on convex body has found a sufficient condition for the solvability of Diophantine inequalities in integers not both zero:
$$ c = {\kappa_p}^p,  \kappa_{p} = \frac{{\Gamma(1 + \frac{2}{p})}^{1/2}}{{\Gamma(1 + \frac{1}{p})}}.$$
 But this result is not optimal, and Minkowski also raised the issue of not improving constant $c$.
For this purpose Minkowski has proposed to use the critical determinant.

Given
any set $\mathcal R\subset {\mathbb R}^n$, a lattice $\Lambda $ is
{\it admissible} for ${\mathcal R}$ (or is $\mathcal R$-{\it
admissible}) if ${\mathcal R}\bigcap\Lambda=\emptyset\mbox{ or
}\{0\}$. The infimum $\Delta(\mathcal R)$ of the determinants (the
determinant of a lattice $\Lambda$ is written $d(\Lambda)$) of all
lattices admissible for $\mathcal R$ is called the {\it critical
determinant} of $\mathcal R$. A lattice $\Lambda $ is {\it
critical} for $\mathcal R$ if $d(\Lambda) = \Delta({\mathcal R})$.

Critical determinant is one of the main 
notion of the geometry
of numbers~\cite{Mi:GZ,Mi:DA,C:GN}. It  has been investigated in the framework of problem of Minkowski in papers by  Mordell~\cite{M:LP}, by Davis~\cite{D:NC},  by Cohn~\cite{Co:MC}, by Watson~\cite{W:MC}, by Malyshev~\cite{Ma:AC,Ma:AC1}
and by Malyshev~\cite{MV:AC,GM:MM,GM:P2,GGM:PM} with colleagues.

\section{Minkowski's conjecture as a problem of Diophantine approximation theory}

Diophantine approximations  connect with critical determinants and with solutions
in integer numbers $x_1, \ldots x_n$  (with some restrictions, for instance  not all  $x_1, \ldots x_n$ are equal to zero)  
of inequalities 
$$  F(x_1, \ldots x_n) < c,$$

or more generally
$$  F(x) < c, \;   x \in \Lambda, x \ne 0.$$

Recall the definitions~\cite{C:GN}.

Let $\mathcal R $ be a set and $\Lambda $ be a  lattice with base 

$ \{a_1, \ldots ,a_n \}$ in ${\bf R}^n.$ A lattice $\Lambda $
is {\it admissible} for body $\mathcal R $ 
($ {\mathcal R}-${\it admissible})
if ${\mathcal D} \bigcap \Lambda = \emptyset $ or $0.$
Let $ d(\Lambda) $ be the determinant of 
$\Lambda.$ The infimum
$\Delta(\mathcal R) $ of determinants of all lattices admissible for
$\mathcal R $ is called {\em the critical determinant} 
of $\mathcal R; $
if there is no $\mathcal R-$admissible lattices then puts
$\Delta(\mathcal R) = \infty. $ A lattice 
$\Lambda $ is {\em critical}
if $ d(\Lambda) = \Delta(\mathcal R).$ \\

 Usually in the
 geometry of numbers the function $F(x)$ is a distance function.
 A real function 
$F(x)$ defined on 
${\bf R}^n$ is {\em distance function} if \\
  
(i) $F(x) \geq 0, x \in {\bf R}^n, F(0) = 0;$ \\
 
(ii) $F(x)$ is continuous;   \\

(iii) $F(x)$ is homogenous: $F(\lambda x) = \lambda F(x), \lambda > 0,  \lambda
\in {\bf R}$. \\

The problem of solving of diophantine inequality $F(x) < c $, 
with a distance function $F$ has the next framework. \\

Let $\overline M$ be the closure of a set $M$ and $\#P$ be the number
of elements of a finite set $P$. An open set 
$ S \subset {\bf R}^n $
is a {\em star body} if $S$ includes the origin of ${\bf R}^n$ and for
any ray $r$ beginning 
in the origin
$ \#(r \cap (\overline M \setminus M)) \leq 1$. If $F(x)$ is a
distance function then the set

$$ M_F = \{ x: F(x) < 1\} $$ 

is a star body.             \\

One of the main particular case of a distance function is the case 
of convex symmetrical function $F(x)$ which with conditions 
(i) - (iii)
 satisfies the additional conditions \\

(iv)$ F(x + y) \leq F(x) + F(y) ;$  \\

(v) $F(-x) = F(x) .$ \\

 The Minkowski's problem can be reformulated as a
 conjecture concerning 
the critical determinant of the region
 $ \mid x \mid^p + \mid y \mid^p  \ \leq 1, \ p > 1.$
 Recall once more that  mentioned mathematical problems are closely 
connected with 
Diophantine Approximation.  \\

 For the given 2-dimension
region $ D_p \subset {\bf R}^2 = (x,y), \ p > 1 $ :

$$ |x|^p + |y|^p < 1 , $$
let $\Delta(D_p) $ 
be the critical determinant of the region.

Let $a \in \Lambda, a \ne 0$ and let
$$  m(F, \Lambda)  = {\inf}_{a} F(a).$$

The Hermite constant of the function $F$ is defined as
$$\gamma(F) = {\sup}_{\Lambda} \frac{ m(F, \Lambda)}{d(\Lambda)^{1/n}}   . $$

\section{Moduli Spaces}

What is moduli? Classically Riemann claimed that $6g - 6$ (real) 
parameters could be for Riemann surface of genus $g > 1$ which would
determine its conformal structure (for elliptic curves, when
$g = 1,$ it is needs one parameter). From algebraic point of view 
we have the following problem: given some kind of variety,
classify the set of all varieties having something in common with
the given one (same numerical invariants of some kind, belonging
to a common algebraic family). For instance, for an elliptic
curve the invariant is the modular invariant of the elliptic
curve.  \\
Let {\bf B} be a class of objects. Let $ S $  be a scheme. 
A family of objects parametrized by
the $ S $ is the  set of objects 
    $ X_{s}: s \in S, X_{s} \in {\bf B} $
 equipped with an additional structure compatible with the structure
  of the base $ S $.   Algebraic moduli 
spaces are defined in the papers by Mumford,  Harris and Morrison
\cite{Mam:EG,HM:MC}. 
\\
 A possibility of the parameterization of all admissible lattices of 
regions $D_p = \{|x|^p + |y|^p < 1\} ,$ under varying $p > 1 ,$ by some
analytical manifolds was  mentioned in the book by Minkowski in 1907 
~\cite{Mi:DA}. 
In 1950 H. Cohn published the paper on
the Minkowski's conjecture \cite{Co:MC}. 
The parameterization
and the corresponding analytic moduli space were
one of the main tools of his approach to the investigation of the 
conjecture. \\ 
Recall some definitions
Let $M$ be an arbitrary  set in ${\bf R}^n$, $O = (0,0) \in {\bf R}^n$.
A lattice
$\Lambda$ is called admissible for $M$, or $M-${\it admissible}, if
it has no points $\neq$ $O$ in the interior of $M.$ It is called
{\it strictly admissible} for
$M$ if it does not contain a point $\neq$ $O$ of $M.$     \\
   The {\it critical determinant} of a set $M$ is the quantity
$\Delta(M)$ given by
$$\Delta(M) = inf\{d(\Lambda): \Lambda \; strictly  \,
admissible \, for M\}  $$
with the understanding that $\Delta(M) = \infty$  if there are no
strictly admissible lattices. The set $M$ is said to be of the finite
or the infinity type according to whether $\Delta(M)$ is finite or
infinite. \\
  The moduli space is defined by the equation
$$ \Delta(p,\sigma) = (\tau + \sigma)(1 + \tau^{p})^{-\frac{1}{p}}
  (1 + \sigma^p)^{-\frac{1}{p}}, \; \; \;  \; (1) $$
in the domain
 $$ {\mathcal M}: \; \infty > p > 1, \; 1 \leq \sigma \leq \sigma_{p} =
 (2^p - 1)^{\frac{1}{p}}, $$
of the $ \{p,\sigma\} $-plane, where $\sigma$ is some real parameter;
$\;$ here $ \tau = \tau(p,\sigma) $ is the function uniquely
determined by the conditions
$$ A^{p} + B^{p} = 1, \; 0 \leq \tau \leq \tau_{p}, $$
where
$$ A = A(p,\sigma) = (1 + \tau^{p})^{-\frac{1}{p}} -
(1 + \sigma^p)^{-\frac{1}{p}},            \;
 B = B(p,\sigma) = \sigma(1 + \sigma^p)^{-\frac{1}{p}}       +
\tau(1 + \tau^{p})^{-\frac{1}{p}}, $$ $\tau_{p}$ is defined by the
equation $ 2(1 - \tau_{p})^{p} = 1 + \tau_{p}^{p}, \; 0 \leq
\tau_{p} < 1. $ \\ 
\begin{definition}
In the notation above, the surface
$$
\Delta-(\tau+\sigma)(1+\tau^{p})^{-{1}/{p}}(1+\sigma^p)^{-{1}/{p}}=0,
$$
in $\mathbf R^3$ with coordinates $(\sigma,p,\Delta)$ we will called
the {\it Minkowski-Cohn moduli space}.
\end{definition}

\section{Minkowski's analytic conjecture}

In considering the question of the minimum value taken by the 
expression $ |x|^p + |y|^p $, with $ p \geq 1 $, 
at points, other
that the origin, of a lattice $\Lambda$ of determinant
$ d(\Lambda)$, Minkowski~\cite{Mi:DA} shows 
that the problem
of determining the maximum value of the minimum for different
lattices may be reduced to that of finding 
the minimum possible
area of a parallelogram with one vertex at the origin and the
three remaining vertices on the curve 
$ |x|^p + |y|^p = 1 $.
The problem with $ p = 1, 2$ and $\infty$ is trivial: in these
cases the minimum areas are $1/2, \: \sqrt{3}/2 $ 
and $1$ 
respectively.
 Let $ D_p \subset {\bf R}^2 = (x,y), \ p > 1 $ be the 2-dimension
region:

$$ |x|^p + |y|^p < 1 . $$
Let $\Delta(D_p) $ 
be the critical determinant of the region.
Recall considerations  of the previous section.
For $p>1$, let 
$$D_p=\{(x,y)\in {\mathbb R}^2\mid|x|^p+|y|^p<1\}.$$
Minkowski~\cite{Mi:DA} raised a question
about critical determinants and critical lattices of regions $D_p$
for varying~$p>1$. Let $\Lambda_{p}^{(0)}$ and $\Lambda_{p}^{(1)}$
be two $D_p$-admissible lattices each of which contains
three pairs of points on the boundary of $D_p$ and with the
property that $(1,0)\in\Lambda_{p}^{(0)},\; (-2^{-1/p},2^{-1/p})
\in \Lambda_{p}^{(1)},$ (under these conditions the lattices are
uniquely defined).
 Using analytic parameterization Cohn~\cite{Co:MC} gives analytic
formulation of Minkowski's conjecture.
 
Let
$$ \Delta(p,\sigma) = (\tau + \sigma)(1 + \tau^{p})^{-\frac{1}{p}}
  (1 + \sigma^p)^{-\frac{1}{p}}, \; \; \;  \; (1) $$

be the function defined 
in the domain
 $$ {\mathcal M}: \; \infty > p > 1, \; 1 \leq \sigma \leq \sigma_{p} =
 (2^p - 1)^{\frac{1}{p}}, $$

of the $ \{p,\sigma\} $ plane, where $\sigma$ 
is some real parameter;
$\;$ here $ \tau = \tau(p,\sigma) $ is the function uniquely
determined by the conditions
 
$$ A^{p} + B^{p} = 1, \; 0 \leq \tau \leq \tau_{p}, $$ 

where
 $$ A = A(p,\sigma) = (1 + \tau^{p})^{-\frac{1}{p}} -
(1 + \sigma^p)^{-\frac{1}{p}}                  $$
 
$$ B = B(p,\sigma) = \sigma(1 + \sigma^p)^{-\frac{1}{p}}       +
\tau(1 + \tau^{p})^{-\frac{1}{p}}, $$
 
$\tau_{p}$ is defined by the equation 

$$ 2(1 - \tau_{p})^{p} = 1 + \tau_{p}^{p}, \; 0 \leq \tau_{p}
\leq 1. $$

 In this case needs to extend the notion of parameter variety to
 parameter manifold. The function $ \Delta(p,\sigma) $ in region
$ {\mathcal M} $  
determines the parameter manifold.


Let ${\Delta^{(1)}_p}  = \Delta(p,1) = 4^{-\frac{1}{p}}\frac{1 +\tau_p }{1 - \tau_p}$,

$ {\Delta^{(0)}_p} = \Delta(p, {\sigma_p}) =  \frac{1}{2}{\sigma}_{p}$. \\

{\bf Minkowski's analytic $(p, \sigma)-$conjecture:} \\

{ \it For any real $p$ with conditions }
$ p > 1, \ p \ne 2, \ 1 < \sigma < \sigma_p $ ,\\

$$ \Delta(p,\sigma) > min({\Delta^{(1)}_p}, {\Delta^{(0)}_p}).$$ \\

In the vicinity of the point $p=1$ and in the vicinity  of the point $(2, \sigma_2)$ the $(p, \tau)$ variant of the Minkowski's analytic conjecture is used.\\

{\bf Minkowski's analytic $(p, \tau)-$conjecture:} \\

Let $\tilde \Delta(p, \tau) = \Delta(p,\sigma), \; \sigma = \sigma(p, \tau): \;  A^{p} + B^{p} = 1$. \\

{ \it For any real $p$ and $\tau$ with conditions }
$ p > 1, \ p \ne 2, \ 0 < \tau < \tau_p $ , \\

$$\tilde \Delta(p,\tau) > min({\Delta^{(1)}_p}, {\Delta^{(0)}_p}).$$ \\

For investigation of properties of function $  \Delta(p,\sigma) $
which are need for proof of Minkowski's conjecture~\cite{Mi:DA,Co:MC} 
we considered
 the value of $\Delta = \Delta(p,\sigma) $ and its
derivatives 
$  \Delta_{\sigma}^{'} \; , \,
\Delta_{\sigma^{2}}^{''} \; , \; \Delta_{p}^{'} \;, \;
 \Delta_{\sigma p}^{''} \; , \; \Delta_{\sigma^{2}p}^{'''} \;$ 
on some subdomains of the domain  $ {\mathcal M} $~\cite{GGM:PM}.

\section{Validated numerics}

  Validated numerics (sometimes called as interval computations)
allow~\cite{M:IA,Sh:IA,AH:IC,BK:AI,SC:WG} \\

(1) rigorous enclosure for roundoff error, truncation error, and
error of data; \\

 (2) computation of rigorous bounds of the ranges of
functions and maps.\\

  A compact closed interval $I = [a,b]$ is the set of real numbers $x$
 such that (s.t.) $a \le x \le b.$
   Interval analysis with this type of intervals uses usually  two sorts
of intervals. Wide intervals are used for representing uncertainty of
the real world or lack of information. Narrow intervals are used for
rounding error bounds. In any of these two cases on each step of an
interval computation we compute the interval $I$ which contains an 
(ideal) solution of our problem. Some examples of implementations of 
the intervals are given in papers~\cite{Gl:I90,Gl:I97}. \\

There are many numerical algorithms for solving mathematical problems.
The majority of these algorithms are iterative, so, since stopping
the algorithms after a certain number of steps, we only get an
approximation $\tilde x$ to the desired solution $x.$ A perfect solution
would if we could estimate the errors of the result not after
the iteration process, but simultaneously with the iteration
process. This is one of the main ideas of interval
analysis~\cite{M:IA,Sh:IA,AH:IC,BK:AI}. \\

\section{Malyshev`s method}

First present each of the expressions
$  \Delta_{\sigma}^{'} \; , \;
\Delta_{\sigma^{2}}^{''} \; , \; \Delta_{p}^{'} \; , \;
 \Delta_{\sigma p}^{''} \; , \; \Delta_{\sigma^{2}p}^{'''} \; $
 in terms of a sum of derivatives of "atoms" 
$ s_{i} = \sigma^{p-i},
\; t_{i} = \tau^{p-i}, \; a_{i} = (1 + \sigma^{p})^{-i-\frac{1}{p}},
\; b_{i} = (1 + \tau^{p})^{-i-\frac{1}{p}}, \; 
A = b_{0} - a_{0},
\; B = \tau b_{0} + \sigma a_{0}, \; \alpha_{i} = A^{p-i}, \;
\beta_{i} = B^{p-i} \; ( i = 0, 1, 2, \ldots).$
  
Then by the implicit function theorem computing
 $ \tau = \tau(p,\sigma) $ by means of the following iteration
 process:

$$ {\tau}_{i + 1} = (1 + {\tau}_{i}^p)^{\frac{1}{p}}
((1 - ((1 + {\tau}_{i}^p)^{-\frac{1}{p}} - (1 + {\sigma}^{p})^
{-\frac{1}{ p}})^{ p})^{\frac{1}{ p}} -  
{ \sigma}(1 + {\sigma}^{p})^{-\frac{1}
{ p}}), $$

For computation of the expression for $\tau_p $
we apply the following iteration:

$$ {(\tau_p)}_{i + 1} = 1 - (2^{-\frac{1}{p}})(1 + {(\tau_p)}_{i}^p)
^{\frac{1}{p}}, \; p > 1, \; {(\tau_p)}_{0} \in [0,0.36].  $$

So we really have represented the function $\Delta$  as the 
function $\Delta(p,\sigma)$ of two variables. The same fact
is true for it's  derivatives.

Recall the construction.
Under given $p$ with increasing $\sigma$ from $ 1$ to ${\sigma}_{p}$ the function $ \tau = \tau(p,\sigma) $  is strictly monotonically decreasing from $\tau_p$  to $0$; $\Delta(p,1) = {\Delta^{(1)}_p} = 4^{-\frac{1}{p}}\frac{1 +\tau_p }{1 - \tau_p}$,  
$\Delta(p, {\sigma_p}) = {\Delta^{(0)}_p} = \frac{1}{2}{\sigma}_{p}$. 

\begin{lemma}
$\frac{\partial \Delta(p,\sigma)}{\partial \sigma}|_{\sigma = 1} = \frac{\partial \Delta(p,\sigma)}{\partial \sigma}|_{\sigma = {\sigma}_{p}} = 0.$
\end{lemma}

A.V. Malyshev noted that technically easier and clearer to work with expressions $l = l(p,\sigma), g = g(p,\sigma), h = h(p,\sigma)$ :

$l = l(p,\sigma) = \Delta(p,\sigma) - \min \{ {\Delta^{(1)}_p} , {\Delta^{(0)}_p} \}$,

$l^{(0)} = l^{(0)}(p,\sigma) = \Delta(p,\sigma) - {\Delta^{(0)}_p} $,

$l^{(1)} = l^{(1)}(p,\sigma) = \Delta(p,\sigma) -  {\Delta^{(1)}_p} $.

$g = g(p,\sigma) =
 -(1 + \sigma^{p})^{1 + \frac{1}{p}} (B^{p-1} - \tau^{p-1}A^{p - 1})\frac{\partial \Delta(p,\sigma)}{\partial \sigma} $.

Then

$sign \frac{\partial \Delta(p,\sigma)}{\partial \sigma} = - sign g(p,\sigma) $.

$h = h(p,\sigma) =  \frac{\partial g(p,\sigma)}{\partial \sigma}$.

These expressions and their derivatives with respect to $\sigma$ and with respect to $p$ are equivalent to expressions
 $  \Delta_{\sigma}^{'} \; , \;
\Delta_{\sigma^{2}}^{''} \; , \; \Delta_{p}^{'} \; , \;
 \Delta_{\sigma p}^{''} \; , \; \Delta_{\sigma^{2}p}^{'''} \; $.

A.V. Malyshev and the author have constructed algebraic expressions for
$\Delta, \Delta_{\sigma}^{'}
 \; , \; \Delta_{\sigma^{2}}^{''} \; , \; \Delta_{p}^{'} \; , \;
 \Delta_{\sigma p}^{''} \; , \; \Delta_{\sigma^{2}p}^{'''} \; $ and respectively 
for $ l(p,\sigma),  g(p,\sigma),  h(p,\sigma)$ and their derivatives with respect to $\sigma$ and with respect to $p$.
Then we at first compute their by fixed point and float point computations.

These calculations have demonstrated the validity of the conjecture.

To prove the theorem for $ p\ge 6$  is enough~\cite{MV:AC} to obtain evaluations 

$h(p,\sigma) < 0$, if $1 \le \sigma \le 1 + \frac{1}{5p}$,

$g(p,\sigma)  < 0$, if $1 + \frac{1}{5p} \le \sigma \le 1 + \frac{1,37}{\sqrt {p}}$,

 $ l(p,\sigma) > 0 $, (i.e. $ \Delta(p,\sigma) > {\Delta^{(1)}_p} $), if  $ 1 + \frac{1,37}{\sqrt {p}} \le \sigma \le {\sigma}_{p} $.

Note that  by Lemma 1 \\
 $ g(p,1) =  g(p,\sigma_p) = 0$.

However, the spread of the approach on the domain $ p < 6 $  met with great difficulties.
To prove this hypothesis in the domain $1 <  p < 6 $ Malyshev proposed, firstly, to extend the class of functions. 

It needs to add derivatives of the functions $ l(p,\sigma) , g(p,\sigma), h(p,\sigma) $    which in themselves are quite complicated. 

These are  derivatives with respect to $\sigma$ and with respect to $p$.

Secondly, Malyshev proposed to construct the interval evaluation of the functions on small intervals, covering the study area.

\subsection{Interval evaluation of functions}

Let $ {\bf X} = ({\bf x}_{1}, \cdots,{\bf x}_{n}) =
([{\underline x}_{1}, {\overline x}_{1}], \cdots,
[{\underline x}_{n}, {\overline x}_{n}] $ 
be the n-dimensional
real interval vector with
 $ {\underline x}_{i} \leq x_{i} \leq {\overline x}_{i} $
 ("rectangle" or "box"). The {\it interval evaluation} of a
 function $ G(x_{1}, \cdots, x_{n}) $ on an interval ${\bf X} $
 is the interval 
$[{\underline G}, {\overline G}] $ such that
for any $ x \in {\bf X}, \; G(x) \in [{\underline G}, {\overline G}]. $
 The interval evaluation is called 
{\it optimal} if
$ {\underline G} = \min G,$ and $ {\overline G} = \max G $
on the interval {\bf X}.   \\
 Let $ D $ be a subdomain of  $ {\mathcal M}. $ Under evaluation in
$ D $ a mentioned function the domain
is covered by rectangles of the form

$$ [{\underline p}, {\overline p}; \;
{\underline \sigma}, {\overline \sigma}]. $$
 
In the case of the formula that expressing 
$  \Delta_{\sigma}^{'} \; , \;
\Delta_{\sigma^{2}}^{''} \; , \; \Delta_{p}^{'} \; , \;
 \Delta_{\sigma p}^{''} \; , \; \Delta_{\sigma^{2}p}^{'''} \; $
 
in terms of a sum of derivatives of "atoms" 
$ s_{i} = \sigma^{p-i},
\; t_{i} = \tau^{p-i}, \; a_{i} = (1 + \sigma^{p})^{-i-\frac{1}{p}},
\; b_{i} = (1 + \tau^{p})^{-i-\frac{1}{p}}, \; A = b_{0} - a_{0},
\; 
B = \tau b_{0} + \sigma a_{0}, \; \alpha_{i} = A^{p-i}, \;
\beta_{i} = B^{p-i} \; ( i = 0, 1, 2, \ldots)$   one applies the
rational interval evaluation to construct formulas
for lower bounds and upper bounds of the functions, which in the
end can be expressed 
in terms of $ {\underline p}, \; {\overline p}, \;
{\underline \sigma}, \; {\overline \sigma}, \; {\underline \tau},
\; {\overline \tau}, \; ; \; $ here 
the bounds $ \; {\underline \tau},
 \; {\overline \tau}, \; $ are obtained with the help of the iteration
process:

$$ {\underline t}_{i + 1} = (1 + {\underline t}_{i}^{\overline p})^
{\frac{1}{\overline p}}((1 - ((1 + {\underline t}_{i}^{\overline p})^
{-\frac{1}{\overline p}}
 - (1 + {\overline \sigma}^{\underline p})^
{-\frac{1}{\underline p}})^{\underline p})^{\frac{1}{\underline p}} -
 {\overline \sigma}(1 + {\overline \sigma}^{\underline p})^{-\frac{1}
{\underline p}}),$$

$$ {\overline t}_{i + 1} = (1 + {\overline t}_{i}^{\underline p})^
{\frac{1}{\underline p}}((1 - ((1 + 
{\overline t}_{i}^{\underline p})^
{-\frac{1}{\underline p}} - (1 + 
{\underline \sigma}^{\overline p})^
{-\frac{1}{\overline p}})^{\overline p})^{\frac{1}{\overline p}} - 
{\underline \sigma}(1 + {\underline \sigma}^{\overline p})^{-\frac{1}
{\overline p}}). $$

$$  \; i = 0,1,\cdots  $$
 
As interval computation is the enclosure method, we have to put:

$$ [{\underline \tau}, \; {\overline \tau}] =
 [{\underline t}_N, \; {\overline t}_N] \bigcap
 [{\underline \tau}_{0}, \; {\overline \tau}_{0}] \; .$$

$ N $ is 
computed on the last step of the iteration. \\

For initial values we may take $: \; [{\underline t}_{0}, \;
{\overline t}_{0}] =
[{\underline \tau}_{0}, \; {\overline \tau}_{0}] = [0,\; 0.36]. $

\subsection{Algorithms and software modules}

Here we give names, input and output of  algorithms and  and software modules 
  for interval evaluation only. All these algorithms and  and software modules  are  implemented, tested  and applied under the computer-assisted proof of Minkowski`s conjecture~\cite{Ma:AC,Ma:AC1,GM:MM,GM:P2,GGM:PM}  . \\

First two algorithms are auxiliary and  described in~\cite{Gl:I90}.

{\bf Algorithm} {\em MonotoneFunction} \\

{\bf Input:} A real function $F(x,y)$ monotonous by $x$ and
by $y.$ \\
Interval $[{\underline x},{\overline x}; {\underline y},
{\overline y}].$ \\

{\bf Output:} The interval evaluation of $F.$    \\

{\bf Algorithm} {\em RationalFunction} \\

{\bf Input:} A rational function $R(x,y).$ 
Interval $[{\underline x},{\overline x}; {\underline y},
{\overline y}].$ \\

{\bf Output:} The interval evaluation of $R.$    \\

Next algorithms and software modules compute functions of Malyshev`s method.

 {\bf Algorithm} {\em TAUPV} \\

{\bf Input:} An implicitly defined function $\tau_p$ from
Section 4. \\
Interval $ [{\underline p}, {\overline p}; \; {\underline \sigma},
{\overline \sigma}].$  \\

{\bf Method:} Iterative interval computation. \\

{\bf Output:} The interval evaluation of $\tau_p$.    \\

{\bf Algorithm} {\em TAUV} \\
{\bf Input:} Implicitly defined function $\tau$ from this
Section. \\

Interval $ [{\underline p}, {\overline p}; \; {\underline \sigma},
{\overline \sigma}].$  \\

{\bf Method:} Described in this Section. \\

{\bf Output:} The interval evaluation of $\tau$.   \\

{\bf Algorithm} {\em L0V} \\

{\bf Input:}  Function $l^0 =  \Delta(p,\sigma) - \Delta_p^{(0)}$.  \\

Interval $ [{\underline p}, {\overline p}; \; {\underline \sigma},
{\overline \sigma}].$  \\

{\bf Method:} Interval computations. \\

{\bf Output:} The interval evaluation of $l^0.$    \\

{\bf Algorithm} {\em L1V} \\

{\bf Input:}  Function $l^1 = \Delta(p,\sigma) - \Delta_p^{(1)}$. \\

Interval $ [{\underline p}, {\overline p}; \; {\underline \sigma},
{\overline \sigma}].$  \\

{\bf Method:} Interval computations. \\

{\bf Output:} The interval evaluation of $l^1.$    \\

{\bf Algorithm} {\em GV} \\

{\bf Input:}  A function $g(p,\sigma)$  which has the same sign
as function $\Delta_{\sigma}^{'}$.\\

Interval $ [{\underline p}, {\overline p}; \; {\underline \sigma},
{\overline \sigma}].$  \\

{\bf Method:} Interval computations. \\

{\bf Output:} The interval evaluation of $g(p,\sigma).$    \\

{\bf Algorithm} {\em HV} \\

{\bf Input:}  A function $h(p,\sigma)$  which is the partial
derivative by $\sigma$ the function $g(p,\sigma)$.\\

$ [{\underline p}, {\overline p}; \; {\underline \sigma},
{\overline \sigma}].$  \\

{\bf Method:} Interval computations. \\
{\bf Output:} The interval evaluation of $h(p,\sigma).$    \\

{\bf Algorithm} {\em DHV} \\

This is the most complicated function and there are several algorithms and software modules for its computation~\cite{GP:89b}.

{\bf Input:}  A function $Dh(p,\sigma)$  which is the partial
derivative by $p$ the function $h(p,\sigma)$.\\

$ [{\underline p}, {\overline p}; \; {\underline \sigma},
{\overline \sigma}].$  \\

{\bf Method:} Interval computations. \\
{\bf Output:} The interval evaluation of $Dh(p,\sigma).$    

\begin{remark}
 $(p, \tau)-$method implemented using algorithms and software modules SIG and SIGV with input parameters respectively      
     $(p, \tau, E1)$     and
   $({\underline p}, {\underline \tau}, {\delta}_p, {\delta}_{\tau}, E1) $.
Here $E1$   is the accuracy of computations of $\sigma$   and $({\underline \sigma}, {\overline \sigma})$;  ${\overline p} = {\underline p} + {\delta}_p, {\overline \tau} = {\underline \tau} + {\delta}_{\tau}$.
\end{remark}

\section{Results}
It is important to note that  Malyshev`s method gives possibility to
prove that a value of the target minimum is an analytic function
but is not a point. Ordinary numerical methods do not allow to
obtain results of the kind.
 
  In notations
\cite{GGM:PM} 
next result have proved:

\begin{theorem}
\cite{GGM:PM}
$$\Delta(D_p) = \left\{
                   \begin{array}{lc}
    \Delta(p,1), \; 1 < p \le 2, \; p \ge p_{0},\\
    \Delta(p,\sigma_p), \;  2 \le p \le p_{0};\\
                     \end{array}
                       \right.
                           $$
here $p_{0}$ is a real number that is defined unique by conditions
$\Delta(p_{0},\sigma_p) = \Delta(p_{0},1),  \;
2,57 \le p_{0}  \le 2,58. $

\end{theorem}

\begin{corollary}
$$ {\kappa_p} = {\Delta(D_p)}^{-\frac{p}{2}}.  $$
\end{corollary}

\section{Strengthened Minkowski`s analytic conjecture}
A.V. Malishev and the author on the base of  some theoretical evidences and results of mentioned computation   have proposed the strengthened Minkowski`s analytic conjecture (MAS)~\cite{GM:MM}.

{\it Strengthened Minkowski's analytic (MAS) conjecture:} \\

{\it For given $ p > 1 $ and increasing $\sigma$ from $1$ to
$\sigma_p$ the function $\Delta(p,\sigma) $ \\

1) increase strictly monotonous if $ 1 < p < 2$ and
$ p \geq p^{(1)} $, \\
 
2) decrease strictly monotonous if $ 2 \leq p \leq p^{(2)} $, \\
 
3) has a unique maximum on the segment $ (1,\sigma_p) $; until the
maximum $\Delta(p,\sigma) $ increase strictly monotonous and 
then decrease 
strictly monotonous if $ p^{(2)} < p < p^{(1)} $; \\

4) constant, if $ p = 2$; \\

here \\

$ p^{(1)} > 2 $ is a root of the equation 
$ \Delta_{\sigma^{2}}^{''}|_{\sigma = \sigma_p } = 0 $; 

$p^{(2)} > 2 $ is a root of the  equation
$ \Delta_{\sigma^{2}}^{''}|_{\sigma = 1} = 0 $ }. \\
 
It  seems that the conjecture (MAS) has not been proved for any parameter except the trivial $ p = 2 $.  \\ 

\section{Conclusions}
 A.V. Malyshev`s approach to Minkowski`s conjecture (in Davis`s  amendment) concerning the critical determinant of the region $|x|^p + |y|^p <1$  for $p > 1$ is proposed  and  A.V. Malyshev`s method of its prove is given.
Applications of the approach and of the method are presented.

\noindent {\it  National Aviation University }

\vskip 2 pt

\noindent {\it e-mail: {\tt glanm@yahoo.com}}


\begin{thebibliography}{99}

\bibitem{Mi:DA} H. Minkowski, {\it Diophantische
Approximationen}, Leipzig: Teubner (1907).

\bibitem{Mi:GZ} H. Minkowski, {\it Geometrie der Zahlen}, Berlin--Leipzig: Teubner (1910).
 
\bibitem{M:LP} L.J. Mordell,  Lattice points in the region $\vert
Ax^4\vert + \vert By^4 \vert \geq 1,$ {\it J. London Math. Soc.}
{\bf 16} , 152--156 (1941).

\bibitem{D:NC} C. Davis,  Note on a conjecture by Minkowski, {\it
J. London Math. Soc.}, {\bf 23}, 172--175 (1948).

\bibitem{C:GN} J. Cassels, {\it An Introduction to the Geometry
of Numbers}, Berlin: Springer-Verlag (1971).

\bibitem{Co:MC} H. Cohn, Minkowski's conjectures on critical lattices
in the metric $\{\vert\xi \vert^p+\vert\eta \vert^p \}^{{1}/{p}},$
{\it Annals of Math.}, {\bf 51}, (2), 734--738 (1950).


\bibitem{W:MC} G. Watson, Minkowski's conjecture on the critical
lattices of the region $|x|^p+|y|^p\leq 1 \;$, (I), (II), {\it
Jour. London Math. Soc.}, {\bf 28}, (3, 4), 305--309, 402--410
(1953).


\bibitem{Ma:AC} A. Malyshev, Application of computers to the proof of a
conjecture of Minkowski's from geometry of numbers. I, {\it Zap.
Nauchn. Semin. LOMI}, {\bf 71},  163--180 (1977).

\bibitem{Ma:AC1} A. Malyshev, Application of computers to the proof of a
conjecture of Minkowski's from geometry of numbers. II, {\it Zap.
Nauchn. Semin. LOMI}, {\bf 82},  29--32 (1979).

\bibitem{MV:AC} A. Malyshev, A. Voronetsky,  The proof of 
 Minkowski's conjecture concerning  the critical determinant of the region $|x|^p+|y|^p<1$ for $p > 6$, {\it Acta arithm.}, 
v. {\bf 71},  447--458 (1975).

\bibitem{GM:MM} N. Glazunov, A. V. Malyshev, 
On Minkowski's critical determinant conjecture,{\it Kibernetika}, No. 5, 10--14 (1985).

\bibitem{GM:P2} N. Glazunov, A. Malyshev. The proof of Minkowski`s conjecture concerning the critical determinant of the 
region $ |x|^p + |y|^p <1 $ near $ p = 2,$(in Russian), {\it Doklady Akad. Nauk  Ukr.SSR} ser.A, 7 .P.9--12 (1986).

\bibitem{GGM:PM} N. Glazunov, A. Golovanov, A. Malyshev,
Proof of Minkowski's hypothesis about the critical determinant of
$|x|^p+|y|^p<1$ domain, {\it Research in Number Theory 9}. Notes
of scientific seminars of LOMI. {\bf 151} Leningrad: Nauka. 40--53
(1986).


\bibitem{Mam:EG} Mumford D.  Towards an Enumerative Geometry of the 
Moduli Space of Curves. {\it Arithmetic and Geometry. Vol. II. Progress
in Math.},  pp.271 - 328, 1983.

\bibitem{HM:MC}
 Harris J., Morrison J., {\it Moduli of curves.}  GTM 187,  Springer, 
Berlin-N.Y, 1998.  


\bibitem{M:IA} R.E. Moore. Interval Analysis. Prentice Hall, Englewood
Cliffs, NJ, 1966.

\bibitem{AH:IC} G. Alefeld, J. Herzberger.
 Introduction to Interval Computations. Academic Press, NY, 1983.

\bibitem{Sh:IA} Yu.I. Shokin. Interval'nij Analiz. Nauka. Seb. District.
  Novosibirsk, 1981 (in Russian).

\bibitem{BK:AI} Applications of Interval Computations (R. Baker Kearfott
and Vladik Kreinovich (Eds.) Kluwer Academic Publlishers. 1996.

\bibitem{SC:WG} Scientific Computing, Validated Numerics, Interval
Methods (Walter Kr{\"a}mer and J{\"u}rgen Wolff von Gudenberg (Eds.))
Kluwer Academic/Plenum Publlishers. NY. 2001. 

\bibitem{GP:89b}  Glazunov N.M. On program package TCHAI (NTAR)
for investigation of functions. Estimates of real functions on
rectangles. (in Russia)
{\it Academy of Sciences of UkSSR. Glushkov Institute of Cybernetics.
Kiev.  VINITI (Moscow)} N3411-B89. 87 p. (1989).

\bibitem{Gl:I90} N.M. Glazunov. Interval arithmetics for evaluation of 
real functions
and its implementation on vector-pipeline computers (in Russian), {\it Issues of
Cybernetics  (Voprosy Kibernetiki). Kernel Software for Supercomputers,}
Moscow.  Acad. of Sci. USSR,  P.91-101 (1990). 

\end{thebibliography}
\end{document}